\documentclass[11pt]{amsart} 
\usepackage[mathscr]{eucal}
\usepackage{amsmath,amsfonts}
\parskip=\smallskipamount

\newtheorem{theorem}{Theorem}[section]

\newtheorem{corollary}[theorem]{Corollary}

\makeatletter
\@addtoreset{equation}{section}
\makeatother

\newcommand{\BB}{{\mathbb B}}
\newcommand{\CC}{{\mathbb C}}

\newcommand{\FF}{{\mathbb F}}

\newcommand{\cD}{{\mathcal D}}
\newcommand{\cE}{{\mathcal E}}

\newcommand{\cG}{{\mathcal G}}
\newcommand{\cH}{{\mathcal H}}
\newcommand{\cK}{{\mathcal K}}

\newcommand{\cP}{{\mathcal P}}

\newcommand{\cY}{{\mathcal Y}}
\newcommand{\cX}{{\mathcal X}}

\newdimen\expt
\def\boxit#1{\setbox0\hbox{$\displaystyle{#1}$}
      \hbox{\lower.4\expt
 \hbox{\lower3\expt\hbox{\lower\dp0
      \hbox{\vbox{\hrule height.4\expt
 \hbox{\vrule width.4\expt\hskip3\expt
      \vbox{\vskip3\expt\box0\vskip2\expt}%
 \hskip3\expt\vrule width.4\expt}\hrule height.4\expt}}}}}}
 
\begin{document}

 \pagestyle{myheadings}
\markboth{ Gelu Popescu}{  Commutator lifting inequalities and interpolation }


\title [  Commutator lifting inequalities and interpolation  ] 
{  Commutator lifting inequalities and interpolation } 
  \author{Gelu Popescu}
\date{October 6, 2003}
\thanks{Research supported in part by an NSF grant}

\address{Department of Mathematics, The University of Texas 
at San Antonio \\ San Antonio, TX 78249, USA}
\email{\tt gpopescu@math.utsa.edu}

\begin{abstract}
In this paper we obtain a multivariable commutator lifting inequality,
 which extends to several variables a recent result of Foia\c s, Frazho, and Kaashoek.
    The inequality yields a  
  multivariable lifting theorem generalizing the noncommutative commutant
  lifting theorem.
  
  This is used to solve  new operator-valued interpolation
   problems of
   Schur-Carath\' eodory, Nevanlinna-Pick, and Sarason type on Fock spaces.
   Some consequences to norm constrained analytic interpolation in the unit
   ball of $\CC^n$ are also considered.

\end{abstract}

\maketitle

\section*{Introduction}
 
In a recent paper \cite{FFK}, 
Foia\c s, Frazho, and Kaashoek
  solved a problem proposed by B. Sz.-Nagy in 1968, extending 
  the commutant lifting theorem \cite{SzF} to
   the case when the underlying operators do not intertwine. 
   Their main result establishes minimal norm
   liftings of certain commutators.
   Our main goal is to obtain a multivariable version of their result.
   
   More precisely,
   let $T:=[T_1,\dots, T_n]$, \ $T_i\in B(\cH)$, be a row contraction, i.e.,
   $$
   T_1T_1^*+\cdots +T_nT_n^*\leq I,$$
    and 
 let 
 $V:=[V_1,\dots, V_n]$, \ $V_i\in B(\cK)$,  be  an isometric lifting of
 $T$ on a Hilbert space $\cK\supseteq \cH$, i.e.,
 $$
 V_i^* V_j=\delta_{ij} I \quad \text{ and  } \quad 
 P_\cH V_i=T_i P_\cH
 $$
  for any \ $i,j=1,\ldots, n$.
 Let $Y:=[Y_1,\dots, Y_n]$, \ $Y_i\in B(\cY)$, be another row contraction
  and  let 
 $W:=[W_1,\dots, W_n]$, \ $W_i\in B(\cX)$, be    an isometric lifting of $Y$ on a 
 Hilbert space  $\cX\supseteq \cY$.
 In Section \ref{commutator}, we prove the following commutator 
 lifting inequality:
 
 {\it
  If
    $A\in B(\cY, \cH)$ is  a 
  contraction, then there is  a contraction $B\in B(\cX, \cK)$ such that
  $ B^*|\cH=A^*$ and 
  \begin{equation*} 
  \|V_iB-BW_i\|\leq \sqrt{2}~ 
  \|[T_1A-AY_1,~\ldots, ~T_nA-AY_n]\|^{1/2}
  \end{equation*}
  for any \ $i=1,\ldots, n$. Moreover,  $\sqrt{2}$  is the best
   possible constant.
   }
 
 Note that in the particular case when $T_iA=AY_i$, \ $i=1,\ldots, n$, 
 the inequality implies the  noncommutative commutant lifting theorem  
   for row contractions \cite{Po-isometric}, 
  \cite{Po-intert} 
   (see \cite{SzF-book}, \cite{FF-book} for the classical 
  case $n=1$). 
When $n=1$, we obtain the  Foia\c s--Frazho--Kaashoek result
\cite{FFK}. In Section \ref{commutator}, we  obtain an improved
   version 
of the above-mentioned inequality (see Theorem \ref{ine1}), 
which has as consequence  a generalization  of 
the noncommutative commutant
  lifting theorem  (see Section \ref{liftings}) and the lifting theorem
   obtained  by  Foia\c s, Frazho, and Kaashoek  
  \cite{FFK1}.
  
 In the last section of this paper, we use our new lifting
  theorem to solve  
   the operator-valued interpolation
   problems of
   Schur-Carath\' eodory (\cite{Sc}, \cite{Ca}), Nevanlinna-Pick \cite{N},
    and Sarason \cite{S} type on Fock spaces.
    
    To give the reader a flavor of our new interpolation results, 
    let us mention (as a particular case) the  scalar
     Nevanlinna-Pick type interpolation problem for 
     $F^2(H_n)$, the  full
    Fock space with $n$ generators.
     
  Let $k,m$ be nonnegative integers and let
  $$
  \BB_n:=\{(\lambda_1,\ldots, \lambda_n)\in \CC^n: 
  \ |\lambda_1|^2+\cdots +|\lambda_n|^2<1\}
  $$
  be the open unit ball of $\CC^n$.
    If $\{z_j\}_{j=1}^m\subset\BB_n$ and $\{w_j\}_{j=1}^m\subset \CC$,
     then there exists  $f\in F^2(H_n)$ such that
     $$
     \|f\|_{\cP_k}\leq 1 \quad \text{ and } \quad f(z_j)=w_j, \quad j=1,\ldots, m,
     $$
     if and only if 
     $$
\left[\frac{ 1 }{ 1-\left< z_j,
z_q\right>}\right]_{j,q=1}^m
\geq 
\left[\frac{ (1- \left< z_j,
z_q\right>^{k+1})w_j \overline{w_q}}{ 1-\left< z_j,
z_q\right>} \right]_{j,q=1}^m.
$$
 We should add that $\|f\|_{\cP_k}$ is defined by
 $$
 \|f\|_{\cP_k}:= \sup  \{\|f\otimes p\|:\   p\in \cP_k, \|p\|\leq 1\},
 $$
 where $\cP_k$ is the set of all polynomials   of degree $\leq k$ 
 in the Fock space $F^2(H_n)$
 and $(F^2(H_n), ~\|\cdot\|_{\cP_k})$  is a Banach space.  Moreover,
 if $f\in F^2(H_n)$ and  $\lim\limits_{k\to \infty} \|f\|_{\cP_k} $ exists, then 
 $f\in F_n^\infty $, the noncommutative analytic Toeplitz algebra introduced in 
 \cite{Po-von} (see also \cite{Po-funct}, \cite{Po-analytic}). In this case
  we have 
 $$
 \|f\|_\infty= \lim_{k\to \infty} \|f\|_{\cP_k}.
 $$
We remark that given $f\in F^2(H_n)$, the {\it evaluation}
 $z \mapsto f(z)$, $z\in \BB_n$,  is a
 holomorphic function in the unit ball of
$\CC^n$ (see Section \ref{interpol}).  
Moreover, Arveson \cite{Arv1} proved that
if $f\in F_n^\infty$, then the map $z \mapsto f(z)$, $z\in \BB_n$,
is a
  multiplier  of the 
reproducing kernel Hilbert space 
with reproducing kernel $K_n: \BB_n\times \BB_n\to \CC$ defined by
 $$
 K_n(z,w):= {\frac {1}
{1-\langle z, w\rangle_{\CC^n}}}, \qquad z,w\in \BB_n.
$$ 
 
 The above-mentioned interpolation problem  is an $F^2(H_n)$ interpolation problem if $k=0$
 and, setting $k\to \infty$, it implies
  the Nevanlinna-Pick 
interpolation problem 
for 
the noncommutative analytic Toeplitz algebra $F^\infty_n $,  which was
 solved by
Arias and the author  
(see \cite{ArPo2}, \cite{Po-interpo}) and, independently, 
 by Davidson and Pitts (see \cite{DP}). 
 Let us mention that, recently, interpolation problems  on 
  the unit  ball   
$\BB_n$  were  also considered in  
\cite{ArPo1}, \cite{AMc2}, \cite{Po-spectral},  \cite{Po-tensor},  
  \cite{BTV}, \cite{Po-meromorphic}, \cite{Po-central}, \cite{Po-nehari}, 
   \cite{BB},  and \cite{EP} .

  In a future paper, we will provide an explicit solution
   (the central interpolant) of our multivariable lifting interpolation
    problem (see Theorem 
   \ref{commutant2}) and show that the maximal entropy principle \cite{FFG}
    is valid 
   in this new setting
   as well as Kaftal-Larson-Weiss  suboptimization type results 
   (\cite{KLW}, \cite{Po-central}) on Fock spaces.
    We will also find explicit solutions
    for the
    operator-valued interpolation
   problems considered in the present paper.

 \section{ Commutator
  lifting  inequalities  
  }\label{commutator}

Let $H_n$ be an $n$-dimensional complex  Hilbert space with orthonormal basis
$e_1$, $e_2$, $\dots,e_n$, where $n\in \{1,2,\dots\}$ or $n=\infty$.
  We consider the full Fock space  of $H_n$ defined by
$$F^2(H_n):=\bigoplus_{k\geq 0} H_n^{\otimes k},$$ 
where $H_n^{\otimes 0}:=\CC 1$ and $H_n^{\otimes k}$ is the (Hilbert)
tensor product of $k$ copies of $H_n$.
Define the left creation 
operators $S_i:F^2(H_n)\to F^2(H_n), \  i=1,\dots, n$,  by
$$
 S_i\psi:=e_i\otimes\psi, \  \psi\in F^2(H_n).
$$
Let $\FF_n^+$ be the unital free semigroup on $n$ generators 
$g_1,\dots,g_n$, and the identity $g_0$.
The length of $\alpha\in\FF_n^+$ is defined by
$|\alpha|:=k$, if $\alpha=g_{i_1}g_{i_2}\cdots g_{i_k}$, and
$|\alpha|:=0$, if $\alpha=g_0$.
We also define
$e_\alpha :=  e_{i_1}\otimes e_{i_2}\otimes \cdots \otimes e_{i_k}$  
 and $e_{g_0}= 1$.
It is  clear that 
$\{e_\alpha:\alpha\in\FF_n^+\}$ is an orthonormal basis of $F^2(H_n)$.
 If $T_1,\dots,T_n\in B(\cH)$  (the algebra of all bounded  
  linear operators on the Hilbert space $ \cH$), define 
$T_\alpha :=  T_{i_1}T_{i_2}\cdots T_{i_k}$,
if $\alpha=g_{i_1}g_{i_2}\cdots g_{i_k}$ and 
$T_{g_0}:=I_\cH$.

Let us recall from \cite{F}, \cite{B},  \cite{Po-models}, 
\cite{Po-isometric}, and  \cite{Po-charact} a few results concerning 
the noncommutative dilation theory for sequences of operators (see
\cite{SzF-book} for the classical case $n=1$).
     A sequence of operators $T:=[T_1,\ldots, T_n]$, $T_i\in B(\cH)$,
        is called
 row contraction if 
 $$
 T_1T_1^*+\cdots +T_nT_n^*\le I_\cH.
 $$
We say that a sequence of isometries   $V:=[V_1,\ldots, V_n]$,\ $V_i\in B(\cK)$,  
  is a  minimal isometric dilation   of $T$  on a Hilbert space 
$\cK\supseteq \cH$ if
 the following properties are satisfied:
\begin{enumerate}
\item[(i)] $V_i^* V_j=0$   ~for all $i\neq j$, \ $i,j\in\{1,\ldots, n\}$;
\item[(ii)]   $V^*_j|\cH=T^*_j$  ~for all  $j=1,\ldots,n$;
\item[(iii)] $\cK=\bigvee\limits_{\alpha\in\FF^+_n}V_\alpha\cH$.
\end{enumerate} 
If $V$ satisfies only the condition (i) and $P_\cH V_i=T_i P_\cH, \ i=1,\dots, n$, 
then $V$ is called isometric lifting of $T$.
The minimal isometric dilation of $T$ is an
 isometric lifting and is uniquely determined up to an isomorphism 
 \cite{Po-isometric}. 

Let us consider a canonical  realization of it on Fock spaces. 
 For convenience of notation, we will sometimes identify the $n$-tuple 
 $T:=[T_1,\ldots, T_n]$ 
 with the row operator  $T:=[T_1~\cdots~T_n]$.
Define  the operator $D_T:\oplus_{j=1}^n\cH\to \oplus_{j=1}^n\cH$ by  setting
$D_T:=(I_{\oplus_{j=1}^n\cH}-T^*T)^{1/2}$, and set
 $\cD:=\overline{D_T(\oplus_{j=1}^n\cH)}$,
 where $\oplus_{j=1}^n \cH$ denotes the direct sum of $n$ copies of $\cH$.
Let $D_i:\cH\to 1\otimes \cD\subset F^2(H_n)\otimes \cD$ be defined by
$$
D_ih:= 1\otimes D_T(\underbrace{0,\ldots,0}_{{i-1}\mbox{ \scriptsize
times}},h,0,\ldots), \quad i=1,\ldots, n.
$$
Consider the Hilbert space $\cK:=\cH\oplus [F^2(H_n)\otimes \cD]$ and 
define $V_i:\cK\to\cK$ by
\begin{equation}\label{dil}
V_i(h\oplus \xi):= T_ih \oplus [D_ih +(S_i\otimes 
I_\cD)\xi]
\end{equation}
for any $h\in \cH, ~\xi\in F^2(H_n)\otimes\cD$.
Note that 
\begin{equation}\label{dilmatr}
V_i=\left[\begin{matrix} T_i& 0\\ 
D_i& S_i\otimes I_\cD
\end{matrix}\right]
\end{equation}
with respect to the decomposition $\cK=\cH\oplus [F^2(H_n)\otimes \cD]$.
In \cite{Po-isometric},  we proved  that  
 $V:=[V_1,\ldots, V_n]$  is the   minimal isometric dilation of $T$.
%
 
 
 The main result of this section is the following  
  lifting inequality
 in several variables.

 \begin{theorem}\label{ine1}
 Let $T:=[T_1,\dots, T_n]$, \ $T_i\in B(\cH)$, be a row contraction and 
 let 
 $V:=[V_1,\dots, V_n]$, \ $V_i\in B(\cK)$,  be  an isometric lifting of
 $T$ on a Hilbert space $\cK\supseteq \cH$.  Let  $\cX_i\subseteq \cX$, 
 \ $i=1,\ldots, n$, be orthogonal subspaces 
  and $R_i\in B(\cX_i,\cX)$ be contractions.
  If  $A\in B(\cX, \cH)$ is  a 
  contraction, then there is  a contraction $B\in B(\cX, \cK)$ such that
  $P_{\cH} B=A$ and 
  \begin{equation}\label{ineq}
  \|V_iBR_i-B|\cX_i\|\leq \sqrt{2}~ 
  \left\|[T_1AR_1-A|\cX_1,~\ldots, ~T_nAR_n-A|\cX_n]\right\|^{1/2}
  \end{equation}
  for any \ $i=1,\ldots, n$. 
  Moreover, $\sqrt{2}$ is the best possible constant.
 \end{theorem}
 \begin{proof}
 Define the operators $X_i\in B(\cX_i, \cH)$ by $X_i:= T_i AR_i$, 
 \ $i=1,\ldots, n$,
 and let $X:=[X_1,\ldots, X_n]$. Since $X$ is a row contraction,  we have
 \begin{equation*}\begin{split}
 \left\|D_X\left(
 \operatornamewithlimits\oplus\limits_{i=1}^n h_i\right)\right\|^2&= 
 \left\|D_T\left(
 \operatornamewithlimits\oplus\limits_{i=1}^n AR_ih_i\right)\right\|^2+
 \left\| 
 \operatornamewithlimits\oplus\limits_{i=1}^n D_A R_ih_i\right\|^2+
 \left\| 
 \operatornamewithlimits\oplus\limits_{i=1}^n  D_{R_i}h_i\right\|^2\\
 &\geq 
 \left\|D_T\left(
 \operatornamewithlimits\oplus\limits_{i=1}^n AR_ih_i\right)\right\|^2+
 \left\| 
 \operatornamewithlimits\oplus\limits_{i=1}^n D_A R_ih_i\right\|^2 
 \end{split}
 \end{equation*} 
for any $h_i\in \cX_i$, $i=1,\ldots, n$. 
 Since $A\in B(\cX, \cH)$ is a contraction and the subspaces $\cX_i$ 
 are orthogonal,  the operator
 $[A|\cX_1,\ldots, A|\cX_n]$ is a contraction acting from the Hilbert space 
 $\operatornamewithlimits\oplus\limits_{i=1}^n \cX_i$ to $\cH$.
Define the operators $M_i\in B(\cX_i, \cH)$ by setting $M_i:=X_i-A|\cX_i$, 
\ $i=1,\ldots, n$, and   denote $M:=[M_1, \ldots, M_n]$.
Since 
\begin{equation}\label{MXA}
M=X-[A|\cX_1,\ldots, A|\cX_n],
\end{equation}
it is clear that $\|M\|\leq 2$. Setting $\gamma:= 2\|M\|$, we have
 $\|M^*M\|\leq \gamma$ and
 it makes sense to define the defect operator
 $D_{M, \gamma}:= (\gamma I-M^*M)^{1/2}\in
  B(\operatornamewithlimits\oplus\limits_{i=1}^n \cX_i)$.
Note also that
\begin{equation}\label{MAMA}
\|[A|\cX_1,\ldots, A|\cX_n]^* M+ M^* [A|\cX_1,\ldots, A|\cX_n]\|\leq \gamma.
\end{equation}
Taking into account relations \eqref{MXA} and \eqref{MAMA}, we obtain
\begin{equation*}
\begin{split}
\left\|D_X\left(
 \operatornamewithlimits\oplus\limits_{i=1}^n h_i\right)\right\|^2&\leq
 \left<
 \left(I-[A|\cX_1,\ldots, A|\cX_n]^*[A|\cX_1,\ldots, A|\cX_n]\right)
 \operatornamewithlimits\oplus\limits_{i=1}^n h_i, 
 \operatornamewithlimits\oplus\limits_{i=1}^n h_i
  \right>\\
 & \qquad + \left< (\gamma I-M^*M) 
 \operatornamewithlimits\oplus\limits_{i=1}^n h_i, 
 \operatornamewithlimits\oplus\limits_{i=1}^n h_i\right>\\
 &= 
 \left\| \operatornamewithlimits\oplus\limits_{i=1}^n h_i\right\|^2-
 \left\| \sum_{i=1}^n Ah_i\right\|^2 +\left\| D_{M,\gamma} \left(
 \operatornamewithlimits\oplus\limits_{i=1}^n h_i\right)\right\|^2\\
 &=
 \left\|D_A\left(\sum_{i=1}^n h_i\right)\right\|^2+
 \left\| D_{M,\gamma} 
 \left(\operatornamewithlimits\oplus\limits_{i=1}^n h_i\right)\right\|^2
\end{split}
\end{equation*}
for any $h_i\in \cX_i$, $i=1,\ldots, n$.
The latter equality is due to the fact that the subspaces $\cX_i$ are orthogonal 
and $A$ is a contraction.  Now, putting  together
 the two   inequalities for $D_X$,   we obtain
$$
\left\| \left[
\begin{matrix}
[D_A |\cX_1,\ldots, D_A|\cX_n]\\
D_{M,\gamma}
\end{matrix}
\right]
\left( \operatornamewithlimits\oplus\limits_{i=1}^n h_i
\right) \right\|\geq  
\left\| \left[
\begin{matrix}
D_ T  \left( \operatornamewithlimits\oplus\limits_{i=1}^n AR_i
\right) \\
\operatornamewithlimits\oplus\limits_{i=1}^n D_AR_i
\end{matrix}
\right]
\left( \operatornamewithlimits\oplus\limits_{i=1}^n h_i
\right) \right\|
$$ 
for any $h_i\in \cX_i$, \ $i=1,\ldots, n$. 
 Hence, and using Douglas factorization theorem, we infer
  that there is a contraction 
 $$
 \left[
 \begin{matrix}
 C& E\\
  Z 
 & F
 \end{matrix}
 \right]:\cD_A\oplus
 \left( \operatornamewithlimits\oplus\limits_{i=1}^n \cX_i
\right)\to 
\cD_T\oplus
 \left( \operatornamewithlimits\oplus\limits_{i=1}^n \cD_A
\right)
 $$
 such that
 $$
 \left[
 \begin{matrix}
 C& E\\
  Z 
 & F
 \end{matrix}
 \right]
 \left[
\begin{matrix}
[D_A |\cX_1,\ldots, D_A|\cX_n]\\
D_{M,\gamma}
\end{matrix}
\right]
=
\left[
\begin{matrix}
D_ T  \left( \operatornamewithlimits\oplus\limits_{i=1}^n AR_i
\right) \\
\operatornamewithlimits\oplus\limits_{i=1}^n D_AR_i
\end{matrix}
\right].
 $$
 Note that the operator $C\in B(\cD_A, \cD_T)$ satisfies the 
 equation
 \begin{equation}\label{l1}
 C[D_A |\cX_1,\ldots, D_A|\cX_n]+ED_{M,\gamma}=D_T
 \left( \operatornamewithlimits\oplus\limits_{i=1}^n AR_i
\right),
 \end{equation}
 and the operator 
 $Z\in B(\cD_A,\operatornamewithlimits\oplus\limits_{i=1}^n \cD_A)$ satisfies
  the equation 
  \begin{equation}\label{l2}
  Z[D_A |\cX_1,\ldots, D_A|\cX_n]+FD_{M,\gamma} =
  \operatornamewithlimits\oplus\limits_{i=1}^n D_AR_i.
  \end{equation}
 The equality \eqref{l1} implies
 \begin{equation}\label{l1'}
 CD_Ah_i+  ED_{M,\gamma, i} h_i = D_i AR_i h_i,\qquad  h_i\in \cX_i,
 \end{equation}
where  $D_{M,\gamma, i}:\cX_i\to \operatornamewithlimits\oplus\limits_{i=1}^n 
\cX_i$ is the $i$ column of the operator matrix  of $ D_{M,\gamma}$.
 Setting 
 $Z=\left[
 \begin{matrix} Z_1\\ \vdots\\Z_n
 \end{matrix}\right]:\cD_A\to
 \operatornamewithlimits\oplus\limits_{i=1}^n \cD_A
 $, 
 the equation \eqref{l2} implies
 $$
 Z_jD_A h_i+P_j FD_{M,\gamma, i}h_i=\delta_{ij} D_AR_ih_i, \qquad h_i\in \cX_i,
 $$
 for any $i,j=1,\ldots, n$, where $P_j$ denotes the orthogonal projection of 
 $\operatornamewithlimits\oplus\limits_{i=1}^n \cD_A$ onto  the $j$ component.
 In particular, if $i=j$, we get
 \begin{equation}\label{Zi}
 Z_iD_A h_i+P_j FD_{M,\gamma, i}h_i= D_AR_ih_i, \qquad h_i\in \cX_i.
 \end{equation}
 Since 
 $\left[
 \begin{matrix} C\\Z_1\\ \vdots\\Z_n
 \end{matrix}\right]:\cD_A\to
 \cD_T\oplus\left(\operatornamewithlimits\oplus\limits_{i=1}^n \cD_A\right)
 $ 
 is a contraction, one can prove that the operator
 $\Lambda:\cD_A\to F^2(H_n)\otimes \cD_T$ defined by
 $$
 \Lambda h:= \sum_{k=0}^\infty \sum_{|\alpha|=k}
  e_\alpha \otimes CZ_{\tilde\alpha} h, \qquad h\in \cD_A,
 $$
 is  also a contraction, where $\tilde\alpha$ stands for the reverse
  of $\alpha=g_{i_1}g_{i_2}\cdots g_{i_k}\in \FF_n^+$, i.e., 
 $\tilde\alpha= g_{i_k}\cdots g_{i_2} g_{i_1}$.
 Indeed, since
 \begin{equation}\label{CZ}
 \|Ch\|^2+ \sum_{i=1}^n \|Z_ih\|^2\leq \|h\|^2, \qquad h\in \cD_A,
 \end{equation}
 we can replace $h$ with $Z_jh$ in \eqref{CZ} and, summing
  up over $j=1,\ldots, n$, we get
  $$\sum_{j=1}^n \|CZ_jh\|^2\leq \sum_{j=1}^n \|Z_jh\|^2-\sum_{i=1}^n 
  \sum_{j=1}^n \|Z_iZ_jh\|^2.
  $$
  Similarly, we obtain
  \begin{equation}\label{CAk}
  \sum_{|\alpha|=k}\|CZ_\alpha h\|^2\leq \sum_{|\alpha|=k}\|Z_\alpha h\|^2-
  \sum_{|\beta|=k+1}\|Z_\beta h\|^2.
  \end{equation}
 Summing up the inequalities \eqref{CAk} for $k=0,1,\ldots, m$, we obtain
 $$
 \sum_{k=0}^m  \sum_{|\alpha|=k}\|CZ_\alpha h\|^2\leq \|h\|^2-
 \sum_{|\beta|=m+1}\|Z_\beta h\|^2\leq \|h\|^2,
 $$
 which proves that $\Lambda$ is a contraction.
 Now, define  the operator 
 $$
 B:\cH\to \cH\oplus [F^2(H_n)\otimes \cD_T]\quad \text{ by } \quad
 B:=\left[\begin{matrix} A\\ \Lambda D_A\end{matrix} \right].
 $$
 We will prove that the contraction $B$ has the required properties. 
  Assume now that $[V_1,\ldots, V_n]$ is the minimal isometric dilation of 
  $[T_1,\ldots, T_n]$.
  Since 
  \begin{equation}\label{Bh}
  Bh=Ah\oplus \sum_{\alpha\in \FF_n^+} e_\alpha\otimes CZ_{\tilde \alpha} D_Ah,
  \qquad h\in \cX,
  \end{equation}
 and taking into account the Fock space realization of
  the minimal isometric dilation of $T$, we obtain
  \begin{equation*}\begin{split}
  V_i BR_ih_i&=
  V_i\left(
  AR_ih_i\oplus \sum_{\alpha\in \FF_n^+} e_\alpha\otimes CZ_{\tilde \alpha}
   D_AR_ih_i
  \right)\\
  &=
  T_i AR_ih_i\oplus\left[
  1\otimes D_iAR_ih_i+\sum_{\alpha\in \FF_n^+} e_{g_i\alpha}\otimes
   CZ_{\tilde \alpha}
   D_AR_ih_i
  \right]
  \end{split}
  \end{equation*}  
 for any $h_i\in \cX_i$,\ $i=1,\ldots, n$.  
 Hence, and using again  relation \eqref{Bh}, we get
  \begin{equation*}\begin{split}
  V_i BR_ih_i-Bh_i=
  (T_i AR_i-A)h_i &\oplus\left[
  1\otimes (D_iAR_i- CD_A)h_i\right.\\
  &\quad+\sum_{\alpha\in \FF_n^+} e_{g_i\alpha}\otimes
   (CZ_{\tilde \alpha}
   D_AR_i-CZ_{\tilde \alpha g_i}D_A)h_i
  ].
  \end{split}
  \end{equation*}  
 Using relations \eqref{MXA}, \eqref{l1'}, and \eqref{Zi}, we obtain
 \begin{equation*}
  V_i BR_ih_i-Bh_i=
  M_ih_i \oplus\left[
  1\otimes ED_{M,\gamma,i}h_i
  \quad+\sum_{\alpha\in \FF_n^+} e_{g_i\alpha}\otimes
   CZ_{\tilde \alpha} P_iFD_{M,\gamma,i}h_i
  \right]
  \end{equation*}  
  for any $h_i\in \cX_i$, \ $i=1,\ldots, n$.
 Hence,
 we deduce
 \begin{equation*} 
 \| V_i BR_ih_i-Bh_i\|^2= 
 \|M_ih_i\|^2 +\|ED_{M,\gamma,i}h_i\|^2
 + \sum_{\alpha\in \FF_n^+}\|CZ_{\tilde \alpha} P_iFD_{M,\gamma,i}h_i\|^2
 \end{equation*}
 for any $h_i\in \cX_i$, \ $i=1,\ldots, n$.
 Since $\Lambda$ and $\left[ \begin{matrix} E\\F\end{matrix} \right]$ are contractions,
 we obtain
 \begin{equation*}\begin{split}
 \| V_i BR_ih_i-Bh_i\|^2 &\leq 
 \|M_ih_i\|^2 +\|ED_{M,\gamma,i}h_i\|^2 +\|P_iFD_{M,\gamma,i}h_i\|^2\\
 &\leq 
 \|M_ih_i\|^2 +\|ED_{M,\gamma,i}h_i\|^2 +\|FD_{M,\gamma,i}h_i\|^2\\
 &\leq \|M_ih_i\|^2 +\|D_{M,\gamma,i}h_i\|^2\\
 &=\|[M_1,\ldots,M_n] k\|^2 +\|D_{M,\gamma}k\|^2\\
 &=\gamma \|h_i\|^2
 \end{split}
 \end{equation*}
 for any $h_i\in \cX_i$, 
 where $k:=\operatornamewithlimits\oplus\limits_{j=1}^n k_j$ with
  $k_i:=h_i$ and $k_j:=0$ if $j\neq i$.
Therefore,
$$
 \| V_i BR_i-B|\cX_i\|\leq \sqrt{\gamma}
 $$
 for any $i=1,\ldots, n$, which proves inequality \eqref{ineq}.
  
  Now, assume that $[V_1,\ldots, V_n]$ is an arbitrary isometric lifting of 
  $[T_1,\ldots, T_n]$.
  The subspace $\cK_0:=\bigvee_{\alpha\in \FF_n^+}V_\alpha \cH$ is 
  reducing under each isometry $V_1,\ldots, V_n$, and 
  $[V_1|\cK_0,\ldots, V_n|\cK_0]$ coincides with the minimal
   isometric dilation of  $[T_1,\ldots, T_n]$. Applying the first part 
   of the proof,
    we find a contraction $B_0\in B(\cX, \cK_0)$ such that $P_\cH B_0=A$  and 
 \begin{equation}\label{K0}
 \| (V_i |\cK_0) B_0R_i-B_0|\cX_i\|\leq \sqrt{\gamma}
 \end{equation}
 for any $i=1,\ldots, n$. Define    $B\in B(\cX,\cK)$ by
  setting  $Bh:=B_0h$,
  ~$h\in \cX$. 
  Now, note that
  $$
  V_i BR_i-B|\cX_i=(V_i |\cK_0) B_0R_i-B_0|\cX_i, \qquad i=1,\ldots, n,
  $$
  and use inequality \eqref{K0}.
  To complete the proof,  notice that the constant $\sqrt{2}$  
  is the best possible  in \eqref{ineq}, since we get equality  for some simple
  examples.
 \end{proof}

 Now, we can prove 
 the commutator
   lifting inequality  announced in the introduction.

 \begin{theorem}\label{ine2}
 Let $T:=[T_1,\dots, T_n]$, \ $T_i\in B(\cH)$, be a row contraction and 
 let 
 $V:=[V_1,\dots, V_n]$, \ $V_i\in B(\cK)$,  be  an isometric lifting of
 $T$ on a Hilbert space $\cK\supseteq \cH$.
 Let $Y:=[Y_1,\dots, Y_n]$, \ $Y_i\in B(\cY)$, be another row contraction
  and  let 
 $W:=[W_1,\dots, W_n]$, \ $W_i\in B(\cX)$, be    an isometric lifting of $Y$ on a 
 Hilbert space  $\cX\supseteq \cY$.
  If
    $A\in B(\cY, \cH)$ is  a 
  contraction, then there is  a contraction $B\in B(\cX, \cK)$ such that
  $ B^*|\cH=A^*$ and 
  \begin{equation}\label{ineq2}
  \|V_iB-BW_i\|\leq \sqrt{2}~ 
  \|[T_1A-AY_1,~\ldots, ~T_nA-AY_n]\|^{1/2}
  \end{equation}
  for any \ $i=1,\ldots, n$.
  Moreover, $\sqrt{2}$ is the best possible constant.
 \end{theorem}
 \begin{proof}
 Define the contraction $\tilde{A}: \cX\to \cH$ by setting
 $$
 \tilde{A}|\cY=A\quad \text{ and } \quad \tilde{A}|\cX\ominus\cY=0.
 $$
 For each $i=1,\ldots, n$, set $\cX_i:= W_i\cX$ and define the operator
  $R_i:\cX_i\to \cX$
   by $R_i:= W_i^*|\cX_i$.
   Since $W_1, \ldots, W_n$ are isometries with orthogonal subspaces, 
   it is clear that the subspaces $\cX_i$ are pairwise orthogonal and $R_i$ are 
   contractions. Applying Theorem \ref{ine1} to the contraction  $\tilde{A}$,
    we find a contraction $B:\cX\to \cK$ such that $P_\cH B=\tilde{A}$ and 
    \begin{equation*}
    \begin{split}
    \|V_iB-BW_i\|&=\|V_iBR_iW_i -BW_i\|\\
    &\leq \|V_iBR_i -B|\cH_i\|\\
    &\leq\sqrt{2} ~
    \|[T_1 \tilde{A}R_1-\tilde{A}|\cX_1,
    \ldots, T_n\tilde{A}R_n-\tilde{A}|\cX_n]\|^{1/2}\\
    &=\sqrt{2}~
    \|[T_1 \tilde{A}-\tilde{A}W_1,
    \ldots, T_n\tilde{A}-\tilde{A}W_n]\|^{1/2}
    \end{split}
    \end{equation*}
    for any \ $i=1,\ldots, n$.
    The latter equality is due to the fact that $[W_1, \ldots, W_n]$ is
     a row
    isometry. Therefore, we  have
    \begin{equation}
    \label{VBW}
    \|V_iB-BW_i\|\leq
    \sqrt{2}
    ~\|[T_1 \tilde{A}-\tilde{A}W_1,
    \ldots, T_n\tilde{A}-\tilde{A}W_n]\|^{1/2}
    \end{equation}
   for any $i=1,\ldots, n$.
   Now, notice that  $ B^*|\cH=A^*$. Moreover, since $[W_1,\ldots, W_n]$
    is an isometric lifting  of $[Y_1,\ldots, Y_n]$ 
    and $\tilde{A}|\cX\ominus\cY=0$, 
    we have
    \begin{equation*}
    \begin{split}
    (T_i\tilde{A}-\tilde{A}W_i)y&=T_i Ay-
    \tilde{A}(P_\cY+P_{\cX\ominus \cY})W_iy\\
    &=T_i Ay-\tilde{A}P_\cY W_iy\\
    &=(T_i A-\tilde{A}Y_i)y
    \end{split}
    \end{equation*}
   for 
   any $y\in \cY$ and $i=1,\ldots, n$.
   On the other hand, we have $(T_i\tilde{A}-\tilde{A}W_i)x=0$
    for any $x\in \cX\ominus \cY$.
    Therefore
    $$
    T_i\tilde{A}-\tilde{A}W_i=[T_i A-AY_i, ~0],\qquad i=1,\ldots, n,
   $$
    with respect to the orthogonal decomposition
    $\cX= \cY\oplus(\cX\ominus \cY)$.
    Using  \eqref{VBW}, we deduce the inequality \eqref{ineq2}.
    The proof is complete.
 \end{proof}
 
 We remark that if one does not require the operator $A$ to be a contraction
 in Theorem \ref{ine2}, then we can find an operator $B$ with the properties
  $\|B\|=\|A\|$,
 $B^*|\cH=A^*$, and 
  \begin{equation*}
  \|V_iB-BW_i\|\leq \sqrt{2}~ \|A\|^{1/2}
  \|[T_1A-AY_1,~\ldots, ~T_nA-AY_n]\|^{1/2}
  \end{equation*}
  for any \ $i=1,\ldots, n$.

 \smallskip
 
  \section{  New lifting  theorems in several variables  
  }\label{liftings}

Let $Y_i\in B(\cY)$ and   $T_i\in B(\cH), \ i=1,\dots, n$,  be operators such that
 $Y:=[Y_1,\dots, Y_n]$  and $T=[T_1,\dots, T_n]$  are row contractions.
Let  $W:=[W_1,\dots, W_n]$  be  an isometric lifting of
 $Y$ on a Hilbert space $\cX\supseteq \cY$, and  $V:=[V_1,\dots, V_n]$  
 be the    an isometric  lifting of $T$ on a Hilbert space $\cK\supseteq 
 \cH$.
Let  $A\in B(\cY, \cH)$ be an operator  
 satisfying $AY_i=T_i A$, \ $i=1,\dots ,n$. An  
  intertwining lifting of $A$ is an operator  
 $B\in B(\cX,\cK)$ satisfying 
$BW_i=V_i B$, \ $i=1,\dots ,n$,  and $P_{\cH} B=AP_\cY$.

 The noncommutative commutant lifting theorem 
  for row contractions
  \cite{Po-isometric}, \cite{Po-intert} 
  (see \cite{SzF}, \cite{FF-book} for the classical case $n=1$) states that if
  $A\in B(\cY, \cH)$
 is an operator satisfying
 $$
 AY_i=T_i A,\qquad i=1,\ldots, n,
 $$
 then there exists an operator $B\in B(\cX,\cK)$ with the following properties:
 \begin{enumerate}
 \item[(i)] $BW_i=V_iB$ ~for any $i=1,\ldots, n$;
 \item[(ii)] $B^*|\cH= A^*$;
 \item[(iii)] $\|B\|= \|A\|$.
 \end{enumerate}
 Note that
 the   noncommutative commutant lifting theorem 
  for row contractions is a consequence of the commutator lifting inequality
  obtained in Theorem \ref{ine2}.

In what follows, we present a  new multivariable
 lifting theorem  which is a simple consequence of    Theorem \ref{ine1}. 
 \begin{theorem}\label{commutant1}
 Let $T:=[T_1,\dots, T_n]$, \ $T_i\in B(\cH)$, be a row contraction and 
 let 
 $V:=[V_1,\dots, V_n]$, \ $V_i\in B(\cK)$,  be  an isometric lifting of
 $T$ on a Hilbert space $\cK\supseteq \cH$.  Let  $\cX_i\subseteq \cX$, 
 \ $i=1,\ldots, n$, be orthogonal subspaces 
  and $R_i\in B(\cX_i,\cX)$ be contractions.
  If  $A\in B(\cX, \cH)$ is   such that 
  $$
  T_iAR_i=A|\cX_i
  $$
  for any  \ $i=1,\ldots, n$, then there is  an operator
     $B\in B(\cX, \cK)$ such that
  $P_{\cH} B=A$, $\|B\|=\|A\|$,  and
  $$
  V_iBR_i=B|\cX_i
  $$ 
  for any \ $i=1,\ldots, n$.
 \end{theorem}

 A very useful equivalent form of Theorem \ref{commutant1} is the following.

 \begin{theorem}\label{commutant2}
 Let $T:=[T_1,\dots, T_n]$, \ $T_i\in B(\cH)$, be a row contraction and 
 let 
 $V:=[V_1,\dots, V_n]$, \ $V_i\in B(\cK)$,  be  an isometric lifting of
 $T$ on a Hilbert space $\cK\supseteq \cH$.  Let $Q_i\in B(\cG_i, \cX)$ 
  be 
 operators with orthogonal ranges and let
 $ C_i \in B(\cG_i, \cX)$ be 
  such that $C_i^* C_i\leq Q_i^* Q_i$ for any $i=1,\ldots, n$. 
  If  $A\in B(\cX, \cH)$ is such that 
  \begin{equation}\label{TAC}
  T_i AC_i= AQ_i 
  \end{equation}
  for any $i=1,\ldots, n$, then there is  an operator
     $B\in B(\cX, \cK)$ such that
  $P_{\cH} B=A$, $\|B\|=\|A\|$ and 
  \begin{equation}\label{VBC}
  V_i BC_i= BQ_i 
  \end{equation}
  for any $i=1,\ldots, n$.
 \end{theorem}
 
 \begin{proof}
 Since $C_i^* C_i\leq Q_i^* Q_i$ for any $i=1,\ldots, n$, there exist
 some 
  contractions
 $R_i:\overline{Q_i \cG_i}\to \cX$ such that $C_i=R_i Q_i$, \ $i=1,\ldots, n$.
 Denote $\cX_i:=\overline{Q_i \cG_i}$ and notice that
  $\cX_i\perp \cX_j$ if $i\neq j$. Note that relation \eqref{TAC} is
   equivalent to $T_iAR_i=A|\cX_i$ for any $i=1,\ldots, n$, 
     and relation \eqref{VBC} is equivalent
   to
   $V_iBR_i=B|\cX_i$ for any $i=1,\ldots, n$. 
   Now, using Theorem \ref{commutant1}, we  can complete the proof.
 \end{proof}
 
 Now let us show that Theorem \ref{commutant2} implies Theorem \ref{commutant1}.
 Indeed, denote $\cG_i:= \cX_i$, $R_i:=C_i$,  and $Q_i:= I_\cX|\cX_i$,
   $i=1,\ldots, n$. Applying Theorem \ref{commutant2}, the implication follows.

 As in the classical case, the general setting of the 
 noncommutative commutant lifting theorem  can be reduced to the case when 
 $Y:=[Y_1,\ldots, Y_n]$ is a row isometry (see \cite{Po-nehari}).
 Notice that Theorem \ref{commutant2} implies the 
  noncommutative commutant lifting theorem. Indeed, is is enough to
  consider $\cG_i:=\cX$, ~$C_i:=I_\cX$, and $Q_i:=Y_i\in B(\cX)$
   for each $i=1,\ldots, n$, where 
  $Y:=[Y_1,\ldots, Y_n]$ is a row isometry.
 
 Applications of Theorem \ref{commutant2} to interpolation on Fock spaces 
 and the unit ball of $\CC^n$  will be considered in the next section.

 \section{Norm constrained interpolation problems on Fock spaces}
 \label{interpol}

   We say that 
 a bounded linear
  operator 
$M\in B(F^2(H_n)\otimes \cK, F^2(H_n)\otimes \cK')$ is  multi-analytic
if 
\begin{equation} \label{a1}
M(S_i\otimes I_\cK)= (S_i\otimes I_{\cK'}) M\quad 
\text{\rm for any }\ i=1,\dots, n.
\end{equation}
Note that $M$ is uniquely determined by the operator
$\theta:\cK\to F^2(H_n)\otimes \cK'$, which is   defined by
 ~$\theta h:=M(1\otimes h)$, \ $h\in \cK$, 
 and
is called the  symbol  of  $M$. We denote $M=M_\theta$. Moreover, 
$M_\theta$ is 
 uniquely determined by the ``coefficients'' 
  $\theta_{(\alpha)}\in B(\cK, \cK')$, which are given by
 \begin{equation}\label{a2}
\left< \theta_{(\tilde\alpha)}h, h'\right>:= \left< \theta h, e_\alpha 
\otimes h'\right>=\left< M_\theta(1\otimes h), e_\alpha 
\otimes h'\right>, 
\end{equation}
where $h\in \cK$, \ $h'\in \cK'$,\ $\alpha\in \FF_n^+$, and 
 $\tilde\alpha$ is the reverse of $\alpha$.
Note that 
$$
\sum\limits_{\alpha \in \FF_n^+} \theta_{(\alpha)}^*   \theta_{(\alpha)}\leq 
\|M_\theta\|^2 I_\cK.
$$
 We can associate with $M_\theta$ a unique formal Fourier expansion 
\begin{equation}\label{four}
M_\theta\sim \sum_{\alpha \in \FF_n^+} R_\alpha \otimes \theta_{(\alpha)},
\end{equation}
where $R_i $, \ $i=1,\ldots, n$, are the right creation operators
on the full Fock space $F^2(H_n)$  
 Since $M_\theta$ acts like its Fourier representation on ``polynomials'', 
  we will identify them for simplicity.
The set of all multi-analytic operators acting  from $F^2(H_n)\otimes \cK$ to
$F^2(H_n)\otimes \cK'$ coincides with $R^\infty_n\bar\otimes B(\cK,\cK')$,
 where $R^\infty_n$ is the weakly closed algebra generated by the right
  creation operators on the full Fock space, and the identity.
 A multi-analytic  operator $M_\theta$ (resp. its symbol $\theta$) is called 
 inner if
 $M_\theta$ is an isometry.
 More about multi-analytic operators on Fock spaces can be found in 
    \cite{Po-charact},
\cite{Po-multi}, \cite{Po-von},  \cite{Po-funct}, and 
 \cite{Po-analytic}.   

 We remark that, in general, if $\theta:\cK\to F^2(H_n)\otimes \cK'$ 
 is a bounded operator
  (which is equivalent to the weak convergence of the series 
 $\sum\limits_{\alpha \in \FF_n^+} \theta_{(\alpha)}^*   \theta_{(\alpha)}$),
  the linear map $M_\theta$ uniquely determined by relations
  \eqref{a2} and \eqref{a1},  
   is not
  a bounded operator. However,  for each $k=0, 1,\ldots$, the  restriction  of 
  $M_\theta$ 
  to $\cP_{k}\otimes \cK$, the set 
  of all polynomias  of degree $\leq k$, is a  bounded operator acting from
  $\cP_{k}\otimes \cK$ to $F^2(H_n)\otimes \cK'$. We define 
  the $\cP_{k}$ norm  of $M_\theta$ by setting
  $$
  \|M_\theta\|_{\cP_{k}}:=\sup\{\|M_\theta q\|: \ q\in \cP_{k}\otimes\cK 
  \text{ and } 
  \|q\|\leq 1\}.
  $$
 It is easy to see that 
  $\|M_\theta\|_{\cP_{k}}\leq \|M_\theta\|_{\cP_{k+1}}$. Note that 
  $M_\theta $ is a multi-analytic operator if and only if
   $\theta\in B(\cK, F^2(H_n)\otimes \cK')$
  and the sequence $\{\|M_\theta\|_{\cP_{k}}\}_{k=0}^\infty$ converges as 
  $k\to\infty$.
  In this case, we have
  $$
  \|M_\theta\|=\lim_{k\to \infty} \|M_\theta\|_{\cP_{k}}.
  $$
  
 For each   $ i=1,\ldots, n$,  
define the operators $C_i$ and $ Q_i$ from $\cP_{k-1}\otimes\cK$ 
to $\cP_{k}\otimes\cK$ by setting
\begin{equation} \label{CQ} 
 C_i:=I_{\cP_{k}\otimes\cK}|\cP_{k-1}\otimes\cK, 
  \quad \text{ and  }\quad 
Q_i:= P_{\cP_{k}\otimes\cK} 
 (S_i\otimes I_\cK) |\cP_{k-1}\otimes\cK,
\end{equation}
where $P_{\cP_{k}\otimes\cK}$ is the orthogonal projection from
 $F^2(H_n)\otimes\cK$ onto $\cP_{k}\otimes\cK$. 
 
 We recall that the invariant subspaces
 under the operators
 $S_1\otimes I_{\cK'}, \ldots, S_n\otimes I_{\cK'}$ ($\cK'$ is a Hilbert space) 
  were characterized
   by the author in \cite{Po-charact}. The next lifting theorem will play an 
 important role in our investigation.

 \begin{theorem}\label{int1}
 Let $\cH\subset F^2(H_n)\otimes \cK'$ be an invariant subspace under each
operator $S_i^*\otimes I_{\cK'}$, $i=1,\ldots, n$, and let
 $A:\cP_k\otimes \cK\to \cH$
 be a bounded operator. 
 Let 
 \begin{equation} \label{Tip}
 T_i:=P_\cH (S_i\otimes I_{\cK'})|\cH, 
 \quad i=1,\ldots, n,
 \end{equation}
  and $C_i, Q_i$ be the operators defined by relation \eqref{CQ}.
 Then there exists an operator $\theta\in B(\cK, F^2(H_n)\otimes \cK')$ 
  such that
 \begin{equation}\label{condi}
 P_\cH M_\theta|\cP_k\otimes \cK= A\quad \text{ and }  
 \quad \|M_\theta\|_{\cP_k}\leq 1
 \end{equation}
 if and only if  ~$\|A\|\leq 1$ and ~$T_iAC_i=AQ_i$~ for any $i=1,\ldots, n$.
 \end{theorem} 
 
 \begin{proof}
 First, note that relation \eqref{Tip} implies
 \begin{equation} \label{PST}
 P_\cH (S_i\otimes I_{\cK'})=T_i P_\cH, \quad i=1,\ldots, n,
 \end{equation}
 which shows that $[S_1\otimes I_{\cK'},\ldots, S_n\otimes I_{\cK'}]$
 is an isometric lifting of $[T_1,\ldots, T_n]$.
 Assume that relation \eqref{condi} holds.  Note that,
  due to the definitions of the operators
 $M_\theta, C_i$ and $Q_i$, we have
 \begin{equation}\label{SMC}
 (S_i\otimes I_{\cK'})(M_\theta|\cP_k\otimes \cK)C_i=
 (M_\theta|\cP_k\otimes \cK)Q_i, \quad i=1,\ldots, n.
 \end{equation} 
 Now, using relations \eqref{PST} and \eqref{SMC}, we obtain
 \begin{equation*}\begin{split}
 T_iAC_i&=T_i (P_\cH M_\theta|\cP_k\otimes \cK) C_i
 = P_\cH (S_i\otimes I_{\cK'})
 (M_\theta|\cP_k\otimes \cK)C_i\\
 &=P_\cH  (M_\theta|\cP_k\otimes \cK)Q_i=AQ_i
 \end{split}
 \end{equation*} 
 for any $i=1,\ldots, n$. It is clear that $\|A\|\leq 1$.
 
 Conversely, assume that $A:\cP_k\otimes \cK\to \cH$ is a contraction
 such that
  ~$T_iAC_i=AQ_i$~ for any $i=1,\ldots, n$.
  According to Theorem \ref{commutant2}, there exists an operator 
  $B:\cP_k\otimes \cK\to
   F^2(H_n)\otimes \cK'$ such that 
   $\|B\|=\|A\|$, $P_\cH B=A$, and
   \begin{equation}\label{SBC}
   (S_i\otimes I_{\cK'})BC_i=BQ_i, \quad i=1,\ldots, n.
   \end{equation}
  Note that if 
   $B:\cP_k\otimes\cK \to F^2(H_n)\otimes \cK'$  is a bounded operator, then 
    there is an operator
  $\theta\in B(\cK,F^2(H_n)\otimes \cK')$ such that 
   $B= M_\theta| \cP_{k}\otimes\cK$
   if and only if relation \eqref{SBC} holds. This completes the proof.
 \end{proof} 
 
 The next result is a Sarason type 
 interpolation theorem \cite{S}  on  Fock spaces, which generalizes the 
 corresponding result for the noncommutative analytic Toeplitz 
 algebra $F_n^\infty$,
   obtained by Arias and the author (see \cite{ArPo2}, \cite{Po-interpo})
    and Davidson and
    Pitts (\cite{DP}).
  
 \begin{theorem}\label{Sarason}
 Let $\varphi\in B(\cK, F^2(H_n)\otimes \cK')$ and let 
 $M_\theta\in R_n^\infty\bar\otimes B(\cE, \cK')$ be an inner 
  multi-analytic operator.
 Then there exists $\psi\in B(\cK, F^2(H_n)\otimes \cE)$ such that
 $$
 \|M_\varphi-M_\theta M_\psi\|_{\cP_k}\leq 1
 $$
 if and only if the operator defined by
 $A:=P_\cH M_\varphi |\cP_k\otimes \cK$ is a contraction, where the
  subspace $\cH$ is defined by
 $$
 \cH:=[F^2(H_n)\otimes \cK']\ominus M_\theta[F^2(H_n)\otimes \cE].
 $$
 \end{theorem}
 \begin{proof}
 First, note that if $f=\varphi-M_\theta \psi$ for some 
 $\psi\in B(\cK, F^2(H_n)\otimes \cE)$, then
 $M_\varphi p-M_f p =M_\theta M_\psi p$ for any polynomial
  $p\in \cP_k\otimes \cK$.
  Hence, 
  $$A=P_\cH M_\varphi |\cP_k\otimes \cK=P_\cH M_f |\cP_k\otimes \cK
  $$ and
  $$
  \|A\|\leq \|M_f\|_{\cP_k}=
  \|M_\varphi - M_\theta M_\psi\|_{\cP_k}.
  $$
  Therefore, we have
  \begin{equation}
  \label{inf}
  \|A\|\leq \inf\{\|M_\varphi - M_\theta M_\psi\|_{\cP_k}:\ 
  \psi\in B(\cK, F^2(H_n)\otimes \cE)\}.
  \end{equation}
  Let us prove that we have equality in \eqref{inf}.
  According to \cite{Po-charact}, the subspace $\cH$ is invariant under each 
  operator 
  $S_i^*\otimes I_{\cK'}$, $i=1,\ldots, n$.
  Let $T_i:=P_\cH (S_i\otimes I_{\cK'})|\cH$, $i=1,\ldots, n$, and note that 
  $[S_1\otimes I_{\cK'},\ldots, S_n\otimes I_{\cK'}]$ is an isometric lifting 
  of 
  $[T_1,\ldots, T_n]$.
   A straighforward calculation shows that 
   \begin{equation*}
   \begin{split}
   T_iAC_i&=T_iP_\cH M_\varphi C_i=P_\cH (S_i\otimes I_{\cK'})M_\varphi C_i\\
   &=P_\cH M_\varphi(S_i\otimes I_{\cK'}) C_i=A(S_i\otimes I_{\cK'}) C_i\\
   &=AQ_i
   \end{split}
   \end{equation*}
   for any $i=1,\ldots, n$.
   Now, using Theorem \ref{int1}, we find $ f\in B(\cK, F^2(H_n)\otimes \cK')$
   such that $A=P_\cH M_f|\cP_k\otimes \cK$ and $\|M_f\|_{\cP_k}=\|A\|$.
  Since $P_\cH(M_\varphi-M_f)|\cP_k\otimes \cK=0$, there exists
  $\psi\in B(\cK, F^2(H_n)\otimes \cE)$ such that $\varphi-f=M_\theta \psi$.
   Hence,
  we deduce
  $$
  \|A\|=\|M_f\|_{\cP_k}=\|M_\varphi -M_\theta\psi\|_{\cP_k},
  $$
  which proves that equality holds in \eqref{inf}.
  This completes the proof.
 \end{proof}

  \begin{corollary}
  Under the hypotheses of Theorem $\ref{Sarason}$, we have
 $$
  \min\{\|M_\varphi-M_\theta M_\psi\|_{\cP_k}:\ 
  \psi\in B(\cK, F^2(H_n)\otimes \cE)\}=\|A\|,
  $$
  where $A:= P_\cH M_\varphi |\cP_k\otimes \cK$.
  \end{corollary}

  Now, we can extend the Schur-Carath\' eodory interpolation result (see 
  \cite{Ca}, \cite{Sc}, and \cite{Po-analytic}) to  Fock spaces.

  \begin{theorem}\label{Cara}
  Let $k, m$ be  nonnegative integers and let 
   $\Theta:=\sum\limits_{|\alpha|\leq m} R_\alpha\otimes \theta_{(\alpha)}$, 
   where
  $\theta_{(\alpha)}\in B(\cK,\cK')$.
  Then there exists an operator $\phi\in B(\cK, F^2(H_n)\otimes \cK')$ such 
  that 
  \begin{equation}\label{the} 
  \theta_{(\alpha)}=\phi_{(\alpha)}\quad \text{if } \ |\alpha|\leq m,
  \quad \text{  and }  \quad \|M_\phi\|_{\cP_k}\leq 1
   \end{equation}
   if and only if
   \begin{equation*}
    \begin{cases} \|P_{\cP_m\otimes\cK'}\Theta|\cP_k\otimes\cK\|
   \leq 1&\quad \text{ if } k\leq m\\
   \|P_{\cP_m\otimes\cK'}\Theta|\cP_m\otimes\cK\|\leq 1&\quad \text{ if } k> m.
   \end{cases}
   \end{equation*}
  \end{theorem}
   \begin{proof}
   Let $\cH:= \cP_k\otimes \cK'$ and note that the subspace $\cH$ is 
   invariant under
    each operator $S_i^*\otimes I_{\cK'}$, $i=1,\ldots, n$, and 
    let 
    $$
    T_i:= P_\cH(S_i\otimes I_{\cK'})|\cH, \quad i=1,\ldots, n.
    $$
    Since $q(S_i\otimes I_{\cK})=(S_i\otimes I_{\cK'})q$ for any $i=1,\ldots, n$,
     a straightforward calculation  shows that 
    \begin{equation}\label{TACA}
    T_iAC_i =A Q_i,\quad i=1,\ldots, n,
    \end{equation}
    where the operator $A:\cP_k\otimes \cK\to \cH$ is defined by
    \begin{equation}\label{PP} A:=
    \begin{cases} P_{\cP_m\otimes\cK'}\Theta|\cP_k\otimes\cK
   &\quad \text{ if } k\leq m\\
   P_{\cP_m\otimes\cK'}\Theta|\cP_m\otimes\cK&\quad \text{ if } k> m.
   \end{cases}
   \end{equation}
   Note that the condition \eqref{the} is equivalent to the existence 
   of an operator \linebreak $\phi\in B(\cK, F^2(H_n)\otimes \cK')$ such that
   \begin{equation*}
   P_\cH M_\phi| \cP_k\otimes \cK=A \quad \text{ and  }\quad
    \|M_\phi| \cP_k\otimes \cK\|\leq 1.
   \end{equation*}
   Now, one can apply Theorem \ref{int1} to complete the proof. 
   \end{proof}  
   
   \begin{corollary} Let $k, m$ be  nonnegative integers and let 
   $\Theta:=\sum\limits_{|\alpha|\leq m} R_\alpha\otimes \theta_{(\alpha)}$, 
   where
  $\theta_{(\alpha)}\in B(\cK,\cK')$.
  Then
  $$
  \min\{ \|M_\phi\|_{\cP_k}: \ \phi\in B(\cK, F^2(H_n)\otimes \cK'), \ 
  \theta_{(\alpha)}=\phi_{(\alpha)}\quad \text{if } \ |\alpha|\leq m\}=\|A\|,
  $$
  where the operator $A$ is defined by relation \eqref{PP}. 
   
   \end{corollary}

 \bigskip

 In what follows, we present 
the left
 tangential Nevanlinna-Pick  interpolation  problem with operatorial 
argument for    $B(\cH, F^2(H_n)\otimes\cK)$.
 
 As in  \cite{Po-models},  the spectral radius associated
 with a  sequence 
$Z:=(Z_1,\ldots, Z_n)$ of operators  $Z_i\in B(\cY)$,   is given by
$$
r(Z):=\lim_{k\to \infty}\left\|\sum_{|\alpha|=k} Z_\alpha 
Z_\alpha^*\right\|^{1/2k}
=\inf_{k\to \infty}\left\|\sum_{|\alpha|=k} Z_\alpha Z_\alpha^*\right\|^{1/2k}.
$$
Note that  if $Z_1Z_1^*+\cdots+ Z_n Z_n^*<rI_\cY$~ with $0<r<1$, 
then $r(Z)<1$.
Any element $\psi \in B(\cH, F^2(H_n)\otimes\cY)$
has a unique 
 representation  
$\psi h:=\sum_{\alpha\in \FF_n^+} e_\alpha\otimes A_{(\alpha)}h$, $h\in \cH$,
for some operators $A_{(\alpha)}\in B(\cH, \cY)$.
Therefore, 
$$
M_\psi\sim \sum_{\alpha\in \FF_n^+} R_{\tilde\alpha}\otimes A_{(\alpha)}
$$
and
$
\|\psi\|^2 =\|\sum_{\alpha\in\FF_n^+}A_{(\alpha)}^* A_{(\alpha)}\|  
$.
If 
$r(Z)<1$,  it makes sense to define {\it the evaluation} of $\psi$
at $(Z_1,\ldots, Z_n)$ by setting
\begin{equation}\label{evaluation}
\psi(Z_1,\ldots, Z_n):=\sum_{k=0}^\infty\sum_{|\alpha|=k}
 Z_{\tilde\alpha} A_{(\alpha)},
\end{equation}
 where $\tilde\alpha$ is the reverse of $\alpha$.
Now, using the fact that the spectral radius of $Z$ is strictly
 less than 1,  one can prove the norm convergence of the series  
 \eqref{evaluation}. 
 Indeed,
it is enough to observe
that
\begin{equation*}\begin{split}
\left\|\sum_{|\alpha|=k} Z_{\tilde\alpha} A_{(\alpha)}\right\|
&\leq 
\left\|\sum_{|\alpha|=k} Z_\alpha Z_{\alpha}^*\right\|^{1/2}
\left\|\sum_{|\alpha|=k} A_{(\alpha)}^* A_{(\alpha)}\right\|^{1/2}\\
&\leq \|f\|~\left\|\sum_{|\alpha|=k} Z_\alpha Z_{\alpha}^*\right\|^{1/2}.
\end{split}
\end{equation*}

Given  $C\in B(\cH, \cY)$, we define 
 $W_{\{Z,C\}}: F^2(H_n)\otimes \cH\to \cY$, the controllability operator
 associated with
$\{Z,C\}$,  by setting
$$
W_{\{Z,C\}} \left(\sum_{\alpha \in\FF_n^+} e_\alpha\otimes h_\alpha\right):= 
\sum_{k=0}^\infty \sum_{|\alpha|=k} Z_{\alpha} C h_\alpha.
$$
Since  $r(Z)<1$, note that $W_{\{Z,C\}}$ is  a well-defined  bounded 
operator. 
 We call the positive operator  $G_{\{Z,C\}}:=W_{\{Z,C\}}W_{\{Z,C\}}^*$ 
   the
 controllability grammian for 
$\{Z,C\}$.
 It is easy to see that 
\begin{equation}\label{gram}
G_{\{Z,C\}}=\sum_{k=0}^\infty\sum_{|\alpha|=k} Z_\alpha CC^*  Z^*_\alpha,
\end{equation}
 where the
series converges in norm. As in the classical case ($n=1$),
 we say that the pair
$\{Z,C\}$ is controllable if its grammian $G_{\{Z,C\}}$ is strictly positive.
We remark that $G_{\{Z,C\}}$ is the unique  positive solution of 
the Lyapunov equation
\begin{equation}\label{L}
X=\sum_{i=1}^n Z_iXZ_i^*+ CC^*.
\end{equation}

For any nonnegative integer $k$, we define 
$W_{\{Z,C, k\}}: \cP_k\otimes \cH\to \cY$, the $k$-controlability operator 
 associated with
$\{Z,C\}$,  by setting
$$
W_{\{Z,C,k\}} \left(\sum_{|\alpha|\leq k } e_\alpha\otimes h_\alpha\right):= 
\sum_{p=0}^k\sum_{|\alpha|=p} Z_{{\alpha}} C h_\alpha.
$$
The corresponding Grammian is 
$G_{\{Z,C,k\}}:=W_{\{Z,C,k\}}W_{\{Z,C,k\}}^*$. 

Let $\cH$,  $\cK$, and  $\cY_i, \ i=1,\ldots, m$,  be Hilbert spaces 
and consider the operators 
\begin{equation} \label{cond}\begin{split}
&B_j:\cK\to \cY_j,\quad  C_j:\cH\to \cY_j,\quad  j=1,\ldots, m\\
& Z_j:=[Z_{j,1},\ldots, Z_{j,n}]:\oplus_{i=1}^n \cY_j\to \cY_j,\quad  j=1,\ldots, m, 
\end{split}
\end{equation}
such that $r(Z_j)<1$ for any $j=1,\ldots, m$.
Given a nonnegative integer $k$,  
the left
 tangential Nevanlinna-Pick  interpolation  problem with operatorial 
argument for    $B(\cH, F^2(H_n)\otimes\cK)$ is 
to find $\phi \in B(\cH, F^2(H_n)\otimes\cK)$ such that $\|M_\phi\|_{\cP_k}\leq 1$
and 
\begin{equation}\label{interp}
[I\otimes B_j) \phi](Z_j)=C_j, \quad  j=1,\ldots, m.
\end{equation}

\begin{theorem}\label{neva} 
Given two nonnegative integer $k,m$, 
the   left tangential Nevanlinna-Pick 
 interpolation  problem with operatorial 
argument 
 and  data 
$Z_j$, $B_j$, and $C_j$, $j=1,\dots, m$, has a solution  in 
     $B(\cH, F^2(H_n)\otimes\cK)$,
if and only if
\begin{equation} \label{np}
\left[\sum_{p=0}^\infty\sum_{|\alpha|=p} Z_{j,\alpha}B_j B_i^* 
 Z_{i,\alpha}^*\right]_{i, j=1}^m\geq\left[\sum_{|\alpha|\leq k}  Z_{j,\alpha}C_j C_i^* 
 Z_{i,\alpha}^*\right]_{i, j=1}^m.
\end{equation}

\end{theorem}
\begin{proof}
Define the following operators
$$
B:=\left[\begin{matrix} B_1\\\vdots\\B_m \end{matrix}\right]:
\cK\to \oplus_{j=1}^m \cY_j,\quad 
C:=\left[\begin{matrix} C_1\\\vdots\\C_m \end{matrix}\right]:
\cH\to \oplus_{j=1}^m \cY_j,
$$
and $Y:=[Y_1,\ldots, Y_n]$, where $Y_i$ is the diagonal operator defined by
$$
Y_i:= \left[\begin{matrix}
Z_{1,i}\\0 \\ \vdots\\ 0\end{matrix}
\begin{matrix}
0\\Z_{2,i}\\\vdots\\ 0\end{matrix}
\begin{matrix}
0\\0\\\vdots\\ Z_{m,i}
\end{matrix}\right]: \oplus_{j=1}^m \cY_j \to \oplus_{j=1}^m \cY_j,
$$
for each $i=1,\ldots, n$.
Since $r(Y)<1$ and   $M_\phi\sim 
\sum\limits_{\alpha\in \FF_n^+} R_{\tilde\alpha}\otimes A_{(\alpha)}$, note that 
\begin{equation}\label{iyc}
[I\otimes B)\phi](Y)=C
\end{equation}
if and only if
 $$
 \sum_{p=0}^\infty\sum_{|\alpha|=p} Z_{j,\tilde\alpha} B_j A_{(\alpha)}=C_j,
 $$
  for any $j=1,\ldots, m$.
 Therefore, relation \eqref{interp} is equivalent to relation \eqref{iyc}.
 On the other hand, a straightforward computation on  the elements of the form
$e_\beta\otimes h$, ~$h\in \cH$, \ $\beta\in \FF_n^+$,  shows that relation
 \eqref{iyc} holds if and only if 
\begin{equation}\label{wpw}
W_{\{Y,B\}}M_\phi|\cP_k\otimes \cH= W_{\{Y,C,k\}}.
\end{equation}
 Now, a simple calculation reveals that
\begin{equation}\label{red}
W_{\{Y,B\}}W^*_{\{Y,B\}}-W_{\{Y,C,k\}}W^*_{\{Y,C,k\}}
=
\sum_{p=0}^\infty\sum_{|\alpha|=p}Y_\alpha BB^* Y_\alpha^* -
 \sum_{p=0}^k\sum_{|\alpha|=p}Y_\alpha CC^* Y_\alpha^*,
\end{equation}
where $W_{\{Y,B\}}$ and $W_{\{Y,C,k\}}$ are the controllability operators
 associated with
$\{Y,B\}$ and $\{Y,C\}$, respectively.
Note that the  inequality
   \eqref{np}  holds
  if and only if    
 \begin{equation}\label{WW*}
 W_{\{Y,B\}}W^*_{\{Y,B\}}-W_{\{Y,C,k\}}W^*_{\{Y,C,k\}}\geq 0.
 \end{equation}

Using the definitions of the controlability operators, we deduce that
\begin{equation}\label{WSi}\begin{split}
W_{\{Y,B\}}(S_i\otimes I_\cK)&= Y_i W_{\{Y,B\}}\\
W_{\{Y,C,k\}}(S_i\otimes I_\cH)| \cP_{k-1}\otimes \cH&= Y_i W_{\{Y,C,k\}}|
 \cP_{k-1}\otimes \cH
\end{split}
\end{equation}
for any $i=1,\ldots, n$.
 Now, it is easy to see that if relation \eqref{wpw} holds and 
 $\| M_\phi|\cP_k\otimes \cH\|\leq 1$, then  the inequality \eqref{WW*} holds. 
 
 Conversely, assume that  inequality \eqref{WW*} holds. Then there exists
  a contraction
  $$
  \Lambda:\overline{\text{\rm range }W_{\{Y,B\}}^*}\to \cP_k\otimes \cH
  $$
   such that 
  $\Lambda W_{\{Y,B\}}^*=W_{\{Y,C,k\}}^*$.
  Since $W_{\{Y,B\}}^* Y_i^*=(S_i^*\otimes I_\cK) W_{\{Y,B\}}^*$  for any  $i=1,\ldots, n$,
  it is clear that the subspace $\cH':=\overline{\text{\rm range }W_{\{Y,B\}}^*}$
  is invariant  under each 
  operator
  $S_i^*\otimes I_\cK$, \ $i=1,\ldots, n$.
   Let 
   \begin{equation}\label{Ti}
   T_i:= 
  P_{\cH'} (S_i\otimes I_\cK)|\cH',\quad i=1,\ldots, n,
  \end{equation}
   and denote
  $A:=\Lambda^*:\cP_k\otimes \cH\to \cH$. Note that 
  \begin{equation}\label{WAW}
  W_{\{Y,B\}}A=W_{\{Y,C,k\}}.
  \end{equation}
  We claim that
  \begin{equation} \label{TACi}
  T_i AC_i=AQ_i,\ i=1,\ldots, n,
  \end{equation}
 where 
 $$C_i:=I_{\cP_{k}\otimes \cH}|\cP_{k-1}\otimes \cH\quad  \text{ and } \quad
 Q_i:= P_{\cP_k\otimes \cH} (S_i\otimes I_\cH)|\cP_{k-1}\otimes \cH
 $$
  for 
 any $i=1,\ldots, n$. Indeed, using relations \eqref{Ti}, \eqref{WAW}, and 
  \eqref{WSi},
   we obtain
   \begin{equation*}
   \begin{split}
   W_{\{Y,B\}} T_i AC_i&= W_{\{Y,B\}}P_{\cH'} (S_i\otimes I_\cK)AC_i=
   W_{\{Y,B\}} (S_i\otimes I_\cK)AC_i\\
   &=Y_i W_{\{Y,B\}} AC_i= Y_i W_{\{Y,C,k\}}C_i\\
   &= W_{\{Y,C,k\}}(S_i\otimes I_\cH)|\cP_{k-1}\otimes \cH=
   W_{\{Y,B\}}AQ_i
   \end{split}
   \end{equation*}
   for 
 any $i=1,\ldots, n$.
 Since $ W_{\{Y,B\}}|\cH'$ is one-to-one, we get relation \eqref{TACi}.
 According to Theorem \ref{int1}, there exists
 $\phi \in B(\cH, F^2(H_n)\otimes\cK)$ such that
  $\|M_\phi\|_{\cP_k}\leq 1$ and 
  \begin{equation}\label{PMP}
  P_{\cH'} M_\phi|\cP_{k}\otimes \cH=A.
  \end{equation}
Using relation \eqref{WAW}, it is easy to see that \eqref{PMP} implies
relation \eqref{wpw}.
 This  completes the proof.
\end{proof}

Now,
assume that $\{Z,B\}$ is controlable, i.e., its grammian $G_{\{Z,B\}}$ 
is strictly positive. It easy to see that the operator $A$ in the proof
 of Theorem \ref{neva} has an explicit formula given by
 \begin{equation}\label{A}
  A:=W_{\{Z,B\}}^*G_{\{Z,B\}}^{-1}G_{\{Z,C,k\}}.
 \end{equation}
\begin{corollary}
Under the conditions of Theorem $\ref{neva}$ and assumimg that
 $\{Z,B\}$ is controlable, we have
$$
 \min\{ \|M_\phi\|_{\cP_k}:\ \phi \in B(\cH, F^2(H_n)\otimes\cK), \ 
[I\otimes B_j) \phi](Z_j)=C_j  
\}=\|A\|,
$$
where the operator $A$ is defined by equation \eqref{A}. 
\end{corollary}

As
 a consequence of Theorem \ref{neva},  we can obtain  
 the following 
  left
 tangential Nevanlinna-Pick   interpolation  problem  
  in the unit ball of $\CC^n$. 
 This extends
 the corresponding results    for the
   noncommutative analytic Toeplitz algebra $F_n^\infty$
   (see \cite{ArPo2}, \cite{Po-interpo}, \cite{DP}, and  \cite{Po-nehari}).

\begin{corollary}
Let $z_j:=(z_{j,1},\ldots, z_{j,n}), \ j=1,\ldots, m$,
 be distinct points in $\BB_n$, 
the open unit ball of $\CC^n$, and let
 $B_j\in B(\cK,\cY_j)$, ~$ C_j\in B(\cH, \cY_j)$, \ $j=1,\ldots, m$,
be bounded operators.
Then, given a nonnegative integer $k$, 
 there exists  an operator $\theta\in B(\cH, F^2(H_n)\otimes \cK)$ such that 
$$
B_j \theta(z_j)=C_j,\quad j=1,\ldots, m,
$$
and 
~$\|M_\theta\|_{\cP_k}\leq 1$ ~if and only if  
$$
\left[\frac{B_j B_q^*}{ 1-\left< z_j,
z_q\right>}\right]_{j,q=1}^m \geq  \left[\frac{ (1- \left< z_j,
z_q\right>^{k+1})C_j C_q^*}{ 1-\left< z_j,
z_q\right>}\right]_{j,q=1}^m.
$$
\end{corollary} 
  \begin{proof}
First, note that for any $j,q=1,\ldots, m$, we have 
$$
\sum_{\alpha\in \FF_n^+} z_{j,\alpha} \bar z_{q,\alpha}=
\frac{ 1}{1-\left< z_j, z_q\right>}.
$$
 In Theorem  \ref{neva}, consider  the particular   case when
$Z_{j,i}:= z_{j,i} I_{\cY_j}, \ j=1,\ldots, m$, and $i=1,\ldots, n$.
A simple computation shows that
$$
G_{\{Z,B\}}=\left[ \sum_{\alpha\in \FF_n^+} z_{j,\alpha}
 \bar z_{q,\alpha}B_j B_q^*\right]_{j, q=1}^m=
 \left[\frac{B_j B_q^*}{1-\left< z_j, z_q\right>}\right]_{j,q=1}^m.
$$

 Hence, we have
\begin{equation*}\begin{split}
G_{\{Z,B\}} -G_{\{Z,C, k\}}&=
 \left[\sum_{p=0}^\infty\sum_{|\alpha|=p} Z_{j,\alpha}B_j B_q^* 
 Z_{q,\alpha}^*\right]_{j,q=1}^m-\left[\sum_{|\alpha|\leq k} 
  Z_{j,\alpha}C_j C_q^* 
 Z_{q,\alpha}^*\right]_{j,q=1}^m\\
&=\left[\frac{B_j B_q^*}{ 1-\left< z_j,
z_q\right>}\right]_{j,q=1}^m -  \left[\frac{ (1- \left< z_j,
z_q\right>^{k+1})C_j C_q^*}{ 1-\left< z_j,
z_q\right>}\right]_{j,q=1}^m..
\end{split}
\end{equation*}
Now, applying Theorem  \ref{neva}, we complete the proof.
  \end{proof}

We remark that the {\it evaluation} $z\mapsto \theta(z)$, $z\in \BB_n$, 
 is an operator-valued holomorphic function  in the unit ball of $\CC^n$.

\end{document}